\documentclass[10pt,a4paper,reqno]{amsart}

\usepackage{avant}
\usepackage{amsmath}
\usepackage{amsfonts}
\usepackage{palatino}
\usepackage{a4wide}
\usepackage{amssymb}
\usepackage{amsthm}

\usepackage{xcolor}
\definecolor{darkgreen}{rgb}{0.8,0.2,0}
\definecolor{citegreen}{rgb}{0,0.6,0}
\definecolor{refred}{rgb}{0.8,0,0}
\usepackage[colorlinks, citecolor=citegreen, linkcolor=refred]{hyperref}

\usepackage[fixamsbugs]{mathtools}
\mathtoolsset{showonlyrefs,showmanualtags}

\usepackage{mathptmx}
\usepackage{mathpazo}
\usepackage{palatino}
\usepackage{a4wide}

\theoremstyle{plain}

\newtheorem{ackn}{Acknowledgements\hspace{-.2em}}

\theoremstyle{definition}

\theoremstyle{remark}

\numberwithin{equation}{section}

\def\NN{{{\mathbb N}}}

\title{A Note on the Lalescu Sequence}

\author[Carlo Mantegazza]{Carlo Mantegazza}
\address[Carlo Mantegazza]{Dipartimento di Matematica e
	Applicazioni ``Renato Caccioppoli'', Universit\`a di Napoli Federico II \& Scuola Superiore Meridionale, Napoli, Italy}

\author[Nicola Pio Melillo]{Nicola Pio Melillo}
\address[Nicola Pio Melillo]{Dipartimento di Matematica e Applicazioni ``Renato Caccioppoli'', Universit\`a di Napoli Federico II\& Scuola Superiore Meridionale, Napoli, Italy}

\date{\today}


\begin{document}

	\begin{abstract}
		The sequence 
		$$
		a_n=\sqrt[n+1]{(n+1)!\,}-\sqrt[n]{n!\,}
		$$
		is called {\em Lalescu sequence}, after the Romanian mathematician Traian Lalescu (1882 -- 1929). We prove that this sequence is monotonically decreasing.
	\end{abstract}

	\maketitle
	
	\section{The Lalescu sequence}
	
	The following sequence 
	$$
	a_n=\sqrt[n+1]{(n+1)!\,}-\sqrt[n]{n!\,}
	$$
	is called {\em Lalescu sequence}, after the Romanian mathematician Traian Lalescu (1882 -- 1929, see~\cite{wiki6}) who proposed it in~\cite{lalescu}, asking about its convergence. Possibly due to its indubitable elegance, on one hand and its not so straightforward analysis, on the other, it attracted various authors, who discussed its properties and generalizations (we underline the evident connection with Stirling's formula and Euler's Gamma function).\\
	The property of decreasing monotonicity was shown by J\'ozsef~S\'andor by a careful analysis of the properties of some functions related to the Gamma function in~\cite{sandor1} (see~\cite[Chapter~5, Section~11]{sandor2} for an English translation), while our proof is a quite complex example of asymptotic analysis and of how one can deal very accurately with the orders of infinitesimals of sequences. Moreover, this line of analysis can actually be applied to a broader class of related sequences.\\
Before starting, let us make clear the connection between the Euler's Gamma function $\Gamma$ and the Lalescu sequence (in particular to its monotonicity). It is well known that $n!=\Gamma(n+1)$, then setting, 
$$
F(x)=\Gamma(x+1)^{1/x}\qquad\text{ and }\qquad A(x)=F(x+1)-F(x),
$$
for every $x>0$, we have $F(n)=\sqrt[n]{n!\,}$ and
$$
a_n=F(n+1)-F(n)=\Gamma(n+2)^{\frac{1}{n+1}}-\Gamma(n+1)^{\frac{1}{n}}=A(n),
$$
Hence, the decreasing monotonicity of $a_n$ would follow from the analogous monotonicity of the function $A$, involving $\Gamma$. Moreover, we can also obtain it if we are able to show that
$$
F(x+2)-F(x+1)<F(x+1)-F(x)\qquad\Longleftrightarrow\qquad \frac{F(x+2)+F(x)}{2}<F(x+1),
$$
for every $x>0$. The last inequality is clearly implied by the strict concavity of $F$, hence of the function 
$$
x\mapsto \Gamma(x+1)^{1/x}.
$$
This is exactly the line followed by S\'andor who proved such concavity in~\cite{sandor1} (\cite[Chapter~5, Section~11]{sandor2}).

	\smallskip
	
	We begin by reviewing some basic facts.
	
	If we suppose that the sequence converges, considering the two sequences $\sqrt[n]{n!}$ and $n$, by means of the Stolz--Cesaro theorem, we have
	$$
	\lim_{n \to \infty}\sqrt[n+1]{(n+1)!\,}-\sqrt[n]{n!\,}=\lim_{n \to \infty}\frac{\sqrt[n]{n!}}{n}=1/e\,,
	$$
	as the last limit is well--known. It actually also follows by the Cesaro--type result (see~\cite{rudin2}, for instance)
	$$
	\lim_{n\to\infty} \frac{x_{n+1}}{x_n}=\ell\qquad\Longrightarrow\qquad\lim_{n\to\infty} \sqrt[n]{x_n}=\ell,
	$$
	considering $x_n=n!/n^n$.\\
	Alternatively, still assuming that the sequence $a_n$ converges to some limit, we have that the sequence given by its arithmetic means
	$$
	\frac{\sum_{k=1}^n a_k}{n}
	$$
	converges to the same limit. So we conclude
	$$
	\lim_{n \to \infty}\sqrt[n+1]{(n+1)!}-\sqrt[n]{n!}=\lim_{n \to \infty}\frac{\sum_{k=1}^n a_k}{n}=\lim_{n \to \infty}\frac{\sqrt[n+1]{(n+1)!}-1}{n}=1/e\,.
	$$
	
	Thus, the tricky part is actually showing that the Lalescu sequence converges. This can be shown by means of Stirling's formula~\cite{Patin,WhiWat,stir1}: we rewrite the sequence as
	$$
	\sqrt[n+1]{(n+1)!\,}-\sqrt[n]{n!\,}=\sqrt[n]{n!}\,\big(e^{\log(n+1)!/(n+1)-\log n!/n}-1\big)
	$$
	and examine the exponent of $e$:
	$$
	\begin{aligned}
		\frac{\log(n+1)!}{n+1}-\frac{\log n!}{n}&=\frac{n\log(n+1)-\log n!}{n(n+1)} \\
		&= \frac{n\log(n+1)+n\log n -n\log n-n\log\sqrt[n]{n!}}{n(n+1)}\\
		&=\frac{n\log(1+1/n)+n\log(n/\sqrt[n]{n!})}{n(n+1)}\\
		&=\frac{\log(1+1/n)+\log(n/\sqrt[n]{n!})}{(n+1)}\,.
	\end{aligned}
	$$
	Then, since $n/\sqrt[n]{n!} \to e$, we have that the exponent is equal to
	$$
	1/(n+1)+o(1/n).
	$$
	If we consider that, by Stirling's formula,
	$$
	\sqrt[n]{n!}\, \approx n/e\,,
	$$
	we get
	$$
	\sqrt[n]{n!}\,\big(e^{\log(n+1)!/(n+1)-\log n!/n}-1\big) \approx \frac{n/e}{n+1}+o(1) \to 1/e\,.
	$$
	
	  An alternative line, without the use of Stirling's formula, goes as follows: we rewrite the sequence as
	$$
	\sqrt[n+1]{(n+1)!\,}-\sqrt[n]{n!\,}=\frac{\sqrt[n]{n!\,}}n\bigg(\frac{\sqrt[n+1]{(n+1)!\,}}{\sqrt[n]{n!\,}}-1\bigg)n
	$$
	and we observe that
	\begin{equation}\label{lim1}
		\bigg(\frac{\sqrt[n+1]{(n+1)!\,}}{\sqrt[n]{n!\,}}\bigg)^n=\bigg(\frac{(n+1)!^n}{n!^{n+1}}\bigg)^{\frac{1}{n+1}}
		=\bigg(\frac{(n+1)^n}{n!}\bigg)^{\frac{1}{n+1}}=\bigg(\frac{n+1}{\sqrt[n]{n!\,}}\bigg)^{\frac{n}{n+1}}\to e\,,
	\end{equation}
as it is well known that $n/\sqrt[n]{n!} \to e$,	which implies	
	\begin{equation}
	\bigg(\frac{n+1}{\sqrt[n]{n!\,}}\bigg)^{\frac{n}{n+1}}=\bigg(\frac{n}{\sqrt[n]{n!\,}}\bigg)^{\frac{n}{n+1}}\bigg(\frac{n+1}{n}\bigg)^{\frac{n}{n+1}}\to e\,.
	\end{equation}
	Then,
	\begin{align*}
		\sqrt[n+1]{(n+1)!\,}-\sqrt[n]{n!\,}=&\,\frac{\sqrt[n]{n!\,}}n\cdot
		\frac{\frac{\sqrt[n+1]{(n+1)!\,}}{\sqrt[n]{n!\,}}-1}{\log\bigg(1+\bigg(\frac{\sqrt[n+1]{(n+1)!\,}}{\sqrt[n]{n!\,}}-1\bigg)\bigg)}\cdot
		\log\bigg(1+\bigg(\frac{\sqrt[n+1]{(n+1)!\,}}{\sqrt[n]{n!\,}}-1\bigg)\bigg)^n\\
		=&\,\frac{\sqrt[n]{n!\,}}n\cdot
		\frac{\frac{\sqrt[n+1]{(n+1)!\,}}{\sqrt[n]{n!\,}}-1}{\log\bigg(1+\bigg(\frac{\sqrt[n+1]{(n+1)!\,}}{\sqrt[n]{n!\,}}-1\bigg)\bigg)}\cdot
		\log\bigg(\frac{\sqrt[n+1]{(n+1)!\,}}{\sqrt[n]{n!\,}}\bigg)^n\to 1/e\,,
	\end{align*}
	for $n\to\infty$. Indeed, the first factor tends to $1/e$, while the second and third go to $1$, by the limit~\eqref{lim1} and 
	\begin{equation}
		\frac{\sqrt[n+1]{(n+1)!\,}}{\sqrt[n]{n!\,}}\to 1\,
	\end{equation}
	which again follows by $n/\sqrt[n]{n!} \to e$.
	
	Another ``natural'' way to show the convergence of the sequence would be to prove that it is bounded and monotone. The boundedness from below is actually easy: the sequence $a_n$ is positive, for every $n\in\NN$. Indeed, when above we expressed the sequence as
	$$
	\sqrt[n+1]{(n+1)!\,}-\sqrt[n]{n!\,}=\sqrt[n]{n!}\,\big(e^{\log(n+1)!/(n+1)-\log n!/n}-1\big)\,,
	$$
	we have seen that the exponent of $e$ is given by
	$$
	\frac{\log(n+1)!}{n+1}-\frac{\log n!}{n}=\frac{\log(1+1/n)+\log(n/\sqrt[n]{n!})}{(n+1)}\,,
	$$
	that is positive, since $1+1/n$ and $n/\sqrt[n]{n!}$ are both greater than $1$, hence the positivity of $a_n$.
	Unfortunately, the monotonicity, that is, the fact that $a_n$ is decreasing (as one could expect), is not present in literature, up to our knowledge and it turns out being absolutely non trivial. 
	
	Our contribution to the study of the Lalescu sequence is then to show such monotonicity, first eventually (from some $n\in\NN$ on, which is clearly sufficient for the convergence) and then fully (for every $n\in\NN$). 
	
	As we will see in the next section, our analysis requires a more refined version of Stirling's formula than the ``standard'' one (a ``higher order'' expansion of $n!$, formula~\eqref{expext}) and some quite precise estimates from above and below on $n!$ (formula~\eqref{robb}). Moreover, to obtain the full monotonicity, some numerical check is also needed in order to deal with the ``small'' values of $n\in\NN$.  
	
	Let us say that we think that what follows can be seen as an interesting (and tough) problem about dealing with orders of infinitesimals by means of Taylor expansions and estimates.

	\section{Decreasing monotonicity}
	
	We set $\ell_n=\sqrt[n]{n!\,}$. Clearly, for every $n \in \mathbb N$, we have $\ell_n>0$.\\
	To see that the sequence 
	$$
	a_n=\sqrt[n+1]{(n+1)!\,}-\sqrt[n]{n!\,}
	$$
	is decreasing, we are going to prove equivalently that
	\begin{equation}\label{1}
		\frac{\ell_{n+2}}{\ell_{n+1}}+\frac{\ell_{n}}{\ell_{n+1}}<2\,. 
	\end{equation}
	Defining
	\begin{equation}
		x_n = \log \frac{\ell_{n+1}}{\ell_n}\,,
	\end{equation}
	the inequality~\eqref{1} can then be written as
	\begin{equation}\label{2}
		\exp(x_{n+1}) + \exp(-x_n) = 2 \exp\Big(\frac{x_{n+1}-x_n}{2}\Big) \cosh \Big(\frac{x_{n+1}+x_n}{2}\Big)<2\,. 
	\end{equation}
	
The ratio $\frac{\ell_{n+1}}{\ell_n}$ is a natural quantity to study, measuring the relative growth between consecutive terms and providing a better understanding of the asymptotic behavior of $\ell_n$, passing to a term to the next one. Moreover, taking the logarithm makes it easier to finely analyze such passage.

	\subsection{Eventual monotonicity}
	
	We are going to use the following ``enhanced'' Stirling's formula (see~\cite{stir1}),
	\begin{align}
		n!=&\,\sqrt{2\pi n\,}\,\Big(\frac{n}{e}\Big)^n\bigg(1+\frac{1}{12n}+\frac{1}{288n^2}+O\Big(\frac{1}{n^3}\Big)\bigg)\\
		=&\,\sqrt{2\pi n\,}\,\Big(\frac{n}{e}\Big)^n\bigg(1+\frac{1}{12n}+\frac{1}{288n^2}+o\Big(\frac{1}{n^2}\Big)\bigg)\,.\label{expext}
	\end{align}
	Then,
	$$
	\log\ell_{n}=\frac{\log(2\pi n)}{2n} +\log n -1 + \frac{1}{n}\log\bigg(1+\frac{1}{12n}+\frac{1}{288n^2}+o\Big(\frac{1}{n^2}\Big)\bigg)\,.
	$$
	Expanding in Taylor series up to $o(1 /n^3)$, we obtain
	\begin{align*}
		\log\ell_{n}=&\,\log n -1 +\frac{\log n}{2n}+\frac{\log(2\pi)}{2n}+\frac{1}{12n^2}- \frac{1}{360n^4}+o\Big(\frac{1} {n^4}\Big)\\
		=&\,\log n -1 +\frac{\log n}{2n}+\frac{\log(2\pi)}{2n}+\frac{1}{12n^2}+o\Big(\frac{1} {n^3}\Big)\,.
	\end{align*}
	Hence, 
	\begin{align*}
		x_n=\log\frac{\ell_{n+1}}{\ell_n}=&\,\log(n+1) -1 +\frac{\log(n+1)}{2(n+1)}+\frac{\log(2\pi)}{2(n+1)}+\frac{1}{12(n+1)^2}+o\Big(\frac{1} {n^3}\Big)\\
		&\,-\log n +1 -\frac{\log n}{2n}-\frac{\log(2\pi)}{2n}-\frac{1}{12n^2}+o\Big(\frac{1} {n^3}\Big)\\
		=&\,\log(1+1/n) +\frac{n\log(1+1/n)-\log n}{2n(n+1)}-\frac{\log(2\pi)}{2n(n+1)}-\frac{2n+1}{12n^2(n+1)^2}+o\Big(\frac{1} {n^3}\Big)\\
		=&\,\log(1+1/n) +\frac{\log(1+1/n)}{2(n+1)}-\frac{\log n}{2n(n+1)}-\frac{\log(2\pi)}{2n(n+1)}-\frac{1}{6n^2(n+1)}+o\Big(\frac{1} {n^3}\Big)\\
		=&\,\frac{1}{n}-\frac{\log n}{2n(n+1)}-\frac{\log(2\pi)}{2n(n+1)}-\frac{1}{2n^2}+\frac{1}{2n(n+1)}\\
		&\,+\frac{1}{3n^3}-\frac{1}{6n^2(n+1)}-\frac{1}{4n^2(n+1)}+o\Big(\frac{1} {n^3}\Big)\\
		=&\,\frac{1}{n}-\frac{\log n}{2n(n+1)}-\frac{\log(2\pi)}{2n(n+1)}-\frac{7}{12n^3}+o\Big(\frac{1} {n^3}\Big)\,.
	\end{align*}
	So we have
	\begin{align*}
		x_{n+1}-x_n=&\,\frac{1}{n+1}-\frac{\log(n+1)}{2(n+1)(n+2)}-\frac{\log(2\pi)}{2(n+1)(n+2)}-\frac{7}{12(n+1)^3}\\
		&\,-\frac{1}{n}+\frac{\log n}{2n(n+1)}+\frac{\log(2\pi)}{2n(n+1)}+\frac{7}{12n^3}+o\Big(\frac{1} {n^3}\Big)\\
		=&\,-\frac{1}{n(n+1)}-\frac{n\log(n+1)-(n+2)\log n}{2n(n+1)(n+2)}+\frac{\log(2\pi)}{n(n+1)(n+2)}+o\Big(\frac{1} {n^3}\Big)\\
		=&\,-\frac{1}{n(n+1)}-\frac{n\log(1+1/n)-2\log n}{2n(n+1)(n+2)}+\frac{\log(2\pi)}{n^3}+o\Big(\frac{1} {n^3}\Big)\\
		=&\,-\frac{1}{n(n+1)}+\frac{\log n}{n^3}-\frac{1}{2n^3}+\frac{\log(2\pi)}{n^3}+o\Big(\frac{1} {n^3}\Big)\\
		=&\,-\frac{1}{n^2}+\frac{\log n}{n^3}+\frac{1}{n^3}-\frac{1}{2n^3}+\frac{\log(2\pi)}{n^3}+o\Big(\frac{1} {n^3}\Big)\\
		=&\,-\frac{1}{n^2}+\frac{\log n}{n^3}+\frac{1+2\log(2\pi)}{2n^3}+o\Big(\frac{1} {n^3}\Big)\,.
	\end{align*}
	For the sum $x_{n+1}+x_n$ we expand in Taylor series up to the term $n^{-2}$.
	\begin{align*}
		x_{n+1}+x_n=&\,\frac{1}{n+1}-\frac{\log(n+1)}{2(n+1)(n+2)}-\frac{\log(2\pi)}{2(n+1)(n+2)}\\
		&\,+\frac{1}{n}-\frac{\log n}{2n(n+1)}-\frac{\log(2\pi)}{2n(n+1)}+o\Big(\frac{1}{n^2}\Big)\\
		=&\,\frac{2n+1}{n(n+1)}-\frac{n\log(n+1)+(n+2)\log n}{2n(n+1)(n+2)}-\frac{\log(2\pi)}{n(n+1)}
		+o\Big(\frac{1}{n^2}\Big)\\
		=&\,\frac{2n+1}{n(n+1)}-\frac{\log[n(n+1)]}{2n(n+1)}-\frac{\log(2\pi)}{n(n+1)}
		+o\Big(\frac{1}{n^2}\Big)\\
		=&\,\frac{2}{n}-\frac{\log[n(n+1)]}{2n(n+1)}-\frac{1+\log(2\pi)}{n^2}
		+o\Big(\frac{1}{n^2}\Big)\,.
	\end{align*}
	It follows, expanding up to $o(1/n^3)$,
	\begin{align*}
		(x_{n+1}+x_n)^2=&\,\bigg(\frac{2}{n}-\frac{\log[n(n+1)]}{2n(n+1)}-\frac{1+\log(2\pi)}{n^2}
		+o\Big(\frac{1}{n^2}\Big)\bigg)^2\\
		=&\,\frac{4}{n^2}-2\frac{\log[n(n+1)]}{n^2(n+1)}-\frac{4(1+\log (2\pi)) }{n^3}+o\Big(\frac{1} {n^3}\Big)\\
		=&\,\frac{4}{n^2}-4\frac{\log n}{n^3}-\frac{4(1+\log (2\pi)) }{n^3}+o\Big(\frac{1} {n^3}\Big)\,.
	\end{align*}
	Developing then in Taylor series formula~\eqref{2}, we obtain
	\begin{align*}
		2 \exp\Big(\frac{x_{n+1}-x_n}{2}&\Big)\,\cosh \Big(\frac{x_{n+1}+x_n}{2}\Big)\\
		=&\,2\bigg(1+\frac{x_{n+1}-x_n}{2}+\dots\bigg)\bigg(1+\frac{1}{2}\Big(\frac{x_{n+1}+x_n}{2}\Big)^2+\dots\bigg)\\
		=&\,2\bigg(1-\frac{1}{2n^2}+\frac{\log n}{2n^3}+\frac{1+2\log(2\pi)}{4n^3}+o\Big(\frac{1}{n^3}\Big)\bigg)\\
		&\,\,\cdot\bigg(1+\frac{1}{8}\bigg(\frac{4}{n^2}-4\frac{\log n}{n^3}-\frac{4(1+\log (2\pi))}{n^3}+o\Big(\frac{1} {n^3}\Big)\bigg)\bigg)\\
		=&\,2\bigg(1-\frac{1}{2n^2}+\frac{\log n}{2n^3}+\frac{1+2\log(2\pi)}{4n^3}+o\Big(\frac{1} {n^3}\Big)\bigg)\\
		&\,\,\cdot\bigg(1+\frac{1}{2n^2}-\frac{\log n}{2n^3}-\frac{1+\log (2\pi)}{2n^3}+o\Big(\frac{1} {n^3}\Big)\bigg)\\
		=&\,2-\frac{1}{2n^3}+o\Big(\frac{1} {n^3}\Big)\,,
	\end{align*}
	which is clearly smaller than $2$, for large $n\in\NN$.
	
	\bigskip
	
	Therefore, the Lalescu sequence is {\em eventually} decreasing, by formula~\eqref{2}.

	\subsection{Full monotonicity}
	
	We will use the following ``standard'' inequalities
	\begin{equation}\label{logdis0}
		x-x^2/2\leqslant\log(1+x)\leqslant x
	\end{equation}
	\begin{equation}\label{logdis1}
		x-x^2/2+x^3/3-x^4/4\leqslant\log(1+x)\leqslant x-x^2/2+x^3/3,
	\end{equation}
	\begin{equation}\label{logdis2}
		1-x\leqslant\frac{1}{1+x}\leqslant 1-x+x^2,
	\end{equation}
	which are valid for $x>0$.\\
	Furthermore,
	\begin{equation}\label{expdis}
		e^x \leqslant \frac{1}{1-x}\qquad\text{ and }\qquad\cosh x\leqslant\frac{1}{1-{x^2}/{2}}\,,
	\end{equation}
	for $x\in(-1,1)$. The second inequality above can be shown by comparing the Taylor series which converge uniformly in the interval $[-1,1]$,
	$$
	\cosh x =\sum_{i=0}^\infty \frac{x^{2n}}{(2n)!}\qquad\text{ and }\qquad\frac{1}{1-{x^2}/{2}}=\sum_{i=0}^\infty \frac{x^{2n}}{2^n}\,,
	$$
	noticing that $2^n\leqslant(2n)!$, for every $n \in \mathbb N$.\\
	From the following estimates due to Robbins~\cite{robbins},
	\begin{equation}\label{robb}
		{\sqrt {2\pi n}}\left({\frac {n}{e}}\right)^{n}e^{\frac {1}{12n+1}}\leqslant n!\leqslant{\sqrt {2\pi n}}\left({\frac {n}{e}}\right)^{n}e^{\frac {1}{12n}}\,,
	\end{equation}
	holding for every $n\in\NN$, it follows that
	$$
	\frac{\log(2\pi n)}{2n} +\log n -1 + \frac{1}{(12n+1)n}\leqslant\log\ell_{n}\leqslant\frac{\log(2\pi n)}{2n} +\log n -1 + \frac{1}{12n^2}\,.
	$$
	Thus,
	\begin{align*}
		x_n=&\,\log\frac{\ell_{n+1}}{\ell_n}\leqslant\log(n+1) -1 +\frac{\log(n+1)}{2(n+1)}+\frac{\log(2\pi)}{2(n+1)}+\frac{1}{12(n+1)^2}\\
		&\,-\log n +1 -\frac{\log n}{2n}-\frac{\log(2\pi)}{2n}-\frac{1}{(12n+1)n}\\
		=&\,\log(1+1/n) +\frac{n\log(1+1/n)-\log n}{2n(n+1)}-\frac{\log(2\pi)}{2n(n+1)}-\frac{23n+12}{12(12n+1)(n+1)^2n}\\
		=&\,\log(1+1/n) +\frac{\log(1+1/n)}{2(n+1)}-\frac{\log n}{2n(n+1)}-\frac{\log(2\pi)}{2n(n+1)}-\frac{23n+12}{12(12n+1)(n+1)^2n}\,.
	\end{align*}
	Then, applying the inequality at the right side of formula~\eqref{logdis1} to $\log(1+1/n)$, we have
	\begin{align*}
		x_n\leqslant&\, \frac{1}{n} -\frac{\log n}{2n(n+1)}-\frac{\log(2\pi)}{2n(n+1)} -\frac{23n+12}{12(12n+1)(n+1)^2n}-\frac{1}{2n^2}+\frac{1}{2n(n+1)}\\
		&\,+ \frac{1}{3n^3}-\frac{1}{4n^2(n+1)}+\frac{1}{6n^3(n+1)}\\
		=&\, \frac{1}{n} -\frac{\log n}{2n(n+1)}-\frac{\log(2\pi)}{2n(n+1)} -\frac{23n+12}{12(12n+1)(n+1)^2n}-\frac{1}{2n^2(n+1)}\\
		&\,+ \frac{1}{12n^3}+\frac{5}{12n^3(n+1)}\\
		=&\, \frac{1}{n} -\frac{\log n}{2n(n+1)}-\frac{\log(2\pi)}{2n(n+1)} -\frac{23n+12}{12(12n+1)(n+1)^2n}- \frac{5}{12n^3}+\frac{11}{12n^3(n+1)}\\
		\leqslant&\, \frac{1}{n} -\frac{\log n}{2n(n+1)}-\frac{\log(2\pi)}{2n(n+1)} -\frac{23n+12}{12(12n+1)(n+1)^2n}- \frac{5}{12n^3}+\frac{11}{12n^4}\,,
	\end{align*}
	where in the last step, we estimated $\frac{1}{12n^3(n+1)}\leqslant\frac{1}{12n^4}$.\\
	By the inequality at the left side of formula~\eqref{logdis2}, we have
	\begin{align*}
		\frac{23}{12(12n+1)(n+1)^2}&=\frac{23}{144n^3\Big(1+\frac{25}{12n}+\frac{7}{6n^2}+\frac{1}{12n^3}\Big)}\\
		&\geqslant \frac{23}{144n^3}\Big(1-\frac{25}{12n}-\frac{7}{6n^2}-\frac{1}{12n^3}\Big) \\
		&\geqslant \frac{23}{144n^3} -\frac{1}{2n^4}\,,
	\end{align*}
	for every $n\geqslant 2$, therefore
	\begin{align}
		x_n\leqslant&\,\frac{1}{n} -\frac{\log n}{2n(n+1)}-\frac{\log(2\pi)}{2n(n+1)}-\frac{23}{144n^3}+\frac{1}{2n^4} -\frac{5}{12n^3}+\frac{11}{12n^4}\\
		\leqslant&\,\frac{1}{n} -\frac{\log n}{2n(n+1)}-\frac{\log(2\pi)}{2n(n+1)}-\frac{83}{144n^3}+\frac{3}{2n^4}\,.\label{x_n1}
	\end{align}
	Furthermore, it is easily seen that, when $n\geqslant 3$, this inequality implies the ``simpler'' inequality
	\begin{equation}
		x_n\leqslant\frac{1}{n} -\frac{\log n}{2n(n+1)}-\frac{\log(2\pi)}{2n(n+1)}\,,\label{x_n10}
	\end{equation}
	which will be useful later.\\
	Similarly, using the left inequalities in formulas~\eqref{logdis0} and~\eqref{logdis1} on $\log(1+1/n)$, we have
	\begin{align*}
		x_n=&\,\log\frac{\ell_{n+1}}{\ell_n}\geqslant\,\log(1+1/n) +\frac{n\log(1+1/n)-\log n}{2n(n+1)}-\frac{\log(2\pi)}{2n(n+1)}-\frac{25n+13}{12(12n+13)(n+1)n^2}\\
		\geqslant&\, \frac{1}{n} -\frac{\log n}{2n(n+1)}-\frac{\log(2\pi)}{2n(n+1)}-\frac{1}{2n^2}+\frac{1}{2n(n+1)}-\frac{25n+13}{12(12n+13)(n+1)n^2}-\frac{1}{4n^2(n+1)}\\
		&+\frac{1}{3n^3}-\frac{1}{4n^4}\\
		=&\,\frac{1}{n} -\frac{\log n}{2n(n+1)}-\frac{\log(2\pi)}{2n(n+1)}-\frac{25n+13}{12(12n+13)(n+1)n^2}-\frac{3}{4n^2(n+1)}+\frac{1}{3n^3}-\frac{1}{4n^4}\\
		\geqslant&\,\frac{1}{n} -\frac{\log n}{2n(n+1)}-\frac{\log(2\pi)}{2n(n+1)}-\frac{25n+13}{12(12n+13)(n+1)n^2}-\frac{3}{4n^3}+\frac{1}{3n^3}-\frac{1}{4n^4}\\
		=&\, \frac{1}{n} -\frac{\log n}{2n(n+1)}-\frac{\log(2\pi)}{2n(n+1})-\frac{25n+13}{12(12n+13)(n+1)n^2}-\frac{5}{12n^3}-\frac{1}{4n^4}
	\end{align*}
	where, in the penultimate step, we estimated $\frac{3}{4n^2(n+1)}\leqslant\frac{3}{4n^3}$.\\
 	Since clearly
	$$
	\frac{25n+13}{12(12n+13)(n+1)n^2}\leqslant \frac{25}{144n^3}+\frac{13}{144n^4}\,,
	$$
	we have
	\begin{align}
		x_n\geqslant&\,\frac{1}{n} -\frac{\log n}{2n(n+1)}-\frac{\log(2\pi)}{2n(n+1)}-\frac{25}{144n^3}-\frac{13}{144n^4}-\frac{5}{12n^3}-\frac{1}{4n^4} \\
		=&\,\frac{1}{n} -\frac{\log n}{2n(n+1)}-\frac{\log(2\pi)}{2n(n+1)}-\frac{85}{144n^3}-\frac{49}{144n^4}\,,\label{x_n2}
	\end{align}
	which implies the ``simpler'' inequality
	\begin{equation}
		x_n\geqslant\frac{1}{n} -\frac{\log n}{2n(n+1)}-\frac{\log(2\pi)}{2n(n+1)}-\frac{1}{n^3}\,,\label{x_n20}
	\end{equation}
	for each $n\geqslant 1$, which we will use later.
	
	Let us then estimate the difference $x_{n+1}-x_n$ from above with the inequalities~\eqref{x_n1} and~\eqref{x_n2}:
	\begin{align*}
		x_{n+1}-x_n\leqslant&\,\frac{1}{n+1}-\frac{\log(n+1)}{2(n+1)(n+2)}-\frac{\log(2\pi)}{2(n+1)(n+2)}-\frac{83}{144(n+1)^3}+\frac{3}{2(n+1)^4}\\
		&\,-\frac{1}{n}+\frac{\log n}{2n(n+1)}+\frac{\log(2\pi)}{2n(n+1)}+\frac{85}{144n^3}+\frac{49}{144n^4}\\
		\leqslant&\,-\frac{1}{n(n+1)}-\frac{n\log(n+1)-(n+2)\log n}{2n(n+1)(n+2)}+\frac{\log(2\pi)}{n(n+1)(n+2)}\\
		&\,+\frac{85(1+3n+3n^2) +2n^3}{144n^3(n+1)^3}+\frac{49}{144n^4}+\frac{3}{2n^4}\\
		\leqslant&\,-\frac{1}{n^2}-\frac{n\log(1+1/n)-2\log n}{2n(n+1)(n+2)}+\frac{\log(2\pi)}{n^3}+\frac{1}{n^3}\\
		&\,+\frac{85(1+3n+3n^2)}{144n^3(n+1)^3}+\frac{2}{144n^3}+\frac{265}{144n^4}\,,
	\end{align*}
	where in the last step we used $-\frac{1}{n(n+1)}=-\frac{1}{n^2}+\frac{1}{n^2(n+1)}\leqslant-\frac{1}{n^2}+\frac{1}{n^3}$.\\
	Since $1+3n+3n^2\leqslant 4n^2$, for $n\geqslant4$, we have
	\begin{align*}
		x_{n+1}-x_n\leqslant&\,\,-\frac{1}{n^2}-\frac{\log(1+1/n)}{2(n+1)(n+2)}+\frac{\log n}{n(n+1)(n+2)}+\frac{\log(2\pi)}{n^3}\\
		&\,+\frac{340}{144n(n+1)^3}+\frac{146}{144n^3}+\frac{265}{144n^4}\\
		\leqslant&\,\,-\frac{1}{n^2}-\frac{\log(1+1/n)}{2(n+1)(n+2)}+\frac{\log n}{n^3}+\frac{\log(2\pi)}{n^3}\\
		&\,+\frac{340}{144n^4}+\frac{146}{144n^3}+\frac{265}{144n^4}\\
		=&\,\,-\frac{1}{n^2}-\frac{\log(1+1/n)}{2(n+1)(n+2)}+\frac{\log n}{n^3}+\frac{\log(2\pi)}{n^3}+\frac{146}{144n^3}+\frac{605}{144n^4}
	\end{align*}
	and applying the left inequality in formula~\eqref{logdis0} to $\log(1+1/n)$, we conclude 
	\begin{align*}
		x_{n+1}-x_n\leqslant&\,-\frac{1}{n^2}-\frac{1}{2n(n+1)(n+2)}+\frac{1}{4n^2(n+1)(n+2)}+\!\frac{\log n}{n^3}+\!\frac{\log(2\pi)}{n^3}+\!\frac{146}{144n^3}+\!\frac{605}{144n^4}\,.\\
		\leqslant&\,-\frac{1}{n^2}-\frac{1}{2n^3}+\frac{3}{2n^4}+\frac{1}{4n^4}+\frac{\log n}{n^3}+\frac{\log(2\pi)}{n^3}+\frac{146}{144n^3}+\frac{605}{144n^4}\,.\\
		\leqslant&\,-\frac{1}{n^2}+\frac{\log n}{n^3}+\frac{74/144+\log(2\pi)}{n^3}+\frac{6}{n^4}\,,
	\end{align*}
	for $n\geqslant4$, where we used $-\frac{1}{2n(n+1)(n+2)}=-\frac{1}{2n^3}+\frac{3n+2}{2n^3(n+1)(n+2)}\leqslant-\frac{1}{2n^3}+\frac{3}{2n^4}$.\\
	On the other hand, using the inequalities~\eqref{x_n10} and~\eqref{x_n20} to estimate $x_{n+1}-x_n$ from below, we have
	\begin{align*}
		x_{n+1}-x_n\geqslant &\,\frac{1}{n+1}-\frac{\log(n+1)}{2(n+1)(n+2)}-\frac{\log(2\pi)}{2(n+1)(n+2)}-\frac{1}{(n+1)^3}\\
		&\,-\frac{1}{n}+\frac{\log n}{2n(n+1)}+\frac{\log(2\pi)}{2n(n+1)}\\
		=&\,-\frac{1}{n(n+1)}-\frac{n\log(n+1)-(n+2)\log n}{2n(n+1)(n+2)}-\frac{1}{(n+1)^3}+\frac{\log(2\pi)}{n(n+1)(n+2)}\\
		=&\,-\frac{1}{n^2}\!+\!\frac{1}{n^2(n+1)}\!-\!\frac{\log(1+1/n)}{2(n+1)(n+2)}\!+\!\frac{\log n}{n(n+1)(n+2)}\!-\!\frac{1}{(n+1)^3}\!+\!\frac{\log(2\pi)}{n(n+1)(n+2)}\\
		\geqslant&\,-\frac{1}{n^2}-\frac{1}{2n(n+1)(n+2)}+\frac{\log n}{n(n+1)(n+2)}+\frac{\log(2\pi)}{n(n+1)(n+2)}\\
		\geqslant&\,-\frac{1}{n^2}\,,
	\end{align*}
	for every $n\in\NN$.\\
	Therefore, for $n\geqslant 4$, we conclude
	\begin{equation}\label{diff-finale}
		-\frac{1}{n^2}\leqslant x_{n+1}-x_n\leqslant-\frac{1}{n^2}+\frac{\log n}{n^3}+\frac{74/144+\log(2\pi)}{n^3}+\frac{6}{n^4}\,.
	\end{equation}
	Furthermore, noticing that the term on the right is certainly negative for $n \geqslant 4$, we have in such case
	\begin{equation}\label{diffquad-finale}
		(x_{n+1}-x_n)^2\leqslant \frac{1}{n^4}\,.
	\end{equation}
	
	Let us now estimate the sum $x_{n+1}+x_n$ using the inequality~\eqref{x_n10}:
	\begin{align*}
		x_{n+1}+x_n\leqslant&\,\frac{1}{n+1}-\frac{\log(n+1)}{2(n+1)(n+2)}-\frac{\log(2\pi)}{2(n+1)(n+2)}+\frac{1}{n}-\frac{\log n}{2n(n+1)}-\frac{\log(2\pi)}{2n(n+1)}\\
		=&\,\frac{2n+1}{n(n+1)}-\frac{n\log(n+1)+(n+2)\log n}{2n(n+1)(n+2)}-\frac{\log(2\pi)}{n(n+2)} \\
		=&\,\frac{2n+1}{n(n+1)}-\frac{\log[n(n+1)]}{2n(n+1)}-\frac{\log(2\pi)}{n(n+2)}+\frac{\log(n+1)}{n(n+1)(n+2)}
	\end{align*}
	and applying inequalities
	\begin{align*}
		&\frac{2n+1}{n(n+1)}=\Big(\frac{2}{n}+\frac{1}{n^2}\Big)\Big(\frac{1}{1+1/n}\Big)\leqslant\frac{2}{n}-\frac{1}{n^2} +\frac{1}{n^3}+\frac{1}{n^4}\leqslant \frac{2}{n}-\frac{1}{n^2} +\frac{2}{n^3}\\
		&\frac{\log[n(n+1)]}{2n(n+1)}\geqslant \frac{\log n}{n(n+1)}\geqslant \frac{\log n}{n^2}-\frac{\log n}{n^3}\\
		&\frac{1}{n(n+1)}\geqslant \frac{1}{n^2}-\frac{1}{n^3}\\
		&\frac{1}{n(n+2)}\geqslant \frac{1}{n^2}-\frac{2}{n^3}
	\end{align*}
	(where in the first one we used the inequality at the right side  of formula~\eqref{logdis2}), we obtain
	$$
	x_{n+1}+x_n\leqslant\frac{2}{n}-\frac{\log n}{n^2} -\frac{1+\log(2\pi)}{n^2}+\frac{\log n}{n^3}+\frac{2+2\log (2\pi)}{n^3}+\frac{\log(n+1)}{n^3}\,,
	$$
	Then, considering the inequality
	\begin{equation}
		\log n+2+2\log (2\pi)+\log(n+1)\leqslant\frac{n}{16}\,,
	\end{equation}
	from which it obviously follows
	\begin{equation}
		\frac{\log n}{n^3}+\frac{2+2\log (2\pi)}{n^3}+\frac{\log(n+1)}{n^3}\leqslant\frac{1}{16n^2}
	\end{equation}
	and which holds for $n \geqslant 271$ (keeping in mind the concavity of the left--hand side and checking numerically), we conclude
	\begin{equation}
		0\leqslant x_{n+1}+x_n\leqslant\frac{2}{n}-\frac{\log n}{n^2} -\frac{15/16+\log(2\pi)}{n^2}\leqslant\frac{2}{n}\,,\label{somma}
	\end{equation}
	for every $n \geqslant 271$.
	
	For $n \geqslant 271$, both $y=\big(\frac{x_{n+1}-x_n}{2}\big)$ and $z=\big(\frac{x_{n+1}+x_n}{2}\big)$ are smaller than $1$ in absolute value, so we can use the inequalities in formula~\eqref{expdis} and evaluate
	$$
	\exp\Big(\frac{x_{n+1}-x_n}{2}\Big)=e^y\leqslant\frac{1}{1-y}=1+y+\frac{y^2}{1-y}\leqslant1+y+y^2\,,
	$$
	since $y\leqslant 0$ and
	$$
	\cosh \Big(\frac{x_{n+1}+x_n}{2}\Big)=\cosh z\leqslant\frac{1}{1-z^2/2}=1+\frac{z^2}{2}+\frac{z^4}{4(1-z^2/2)}\leqslant 1+\frac{z^2}{2}+\frac{z^4}{2}\,,
	$$
	since $z^2\leqslant 1$. Therefore, for the inequalities~\eqref{diff-finale},~\eqref{diffquad-finale} and~\eqref{somma}, we get
	\begin{align*}
		2 \exp\Big(&\frac{x_{n+1}-x_n}{2}\Big)\,\cosh \Big(\frac{x_{n+1}+x_n}{2}\Big)\leqslant\big(1+y+y^2\big)\,\Big(1+\frac{z^2}{2}+\frac{z^4}{2}\Big)\\
		=&\, 2\bigg(1+\frac{x_{n+1}-x_n}{2}+\Big(\frac{x_{n+1}-x_n}{2}\Big)^2\bigg)\bigg(1+\frac{1}{2} \Big(\frac{x_{n+1}+x_n}{2}\Big)^2+\frac{1}{2} \Big(\frac{x_{n+1}+x_n}{2}\Big)^4\bigg)\\
		\leqslant&\, 2\bigg(1-\frac{1}{2n^2}+\frac{\log n}{2n^3}+\frac{74/144+\log(2\pi)}{2n^3}+\frac{3}{n^4}+\frac{1}{4n^4}\bigg)\\
		&\,\,\cdot\bigg(1+\frac{1}{2}\Big(\frac{1}{n}-\frac{\log n}{2n^2} -\frac{15/16+\log(2\pi)}{2n^2}\Big)^2+\frac{1}{2n^4}\bigg)\\
		=&\,2\bigg(1-\frac{1}{2n^2}+\frac{\log n}{2n^3}+\frac{74/144+\log(2\pi)}{2n^3}+\frac{13}{4n^4}\bigg)\\
		&\,\,\cdot\bigg(1+\frac{1}{2n^2}-\frac{15/16+\log(2\pi)+\log n}{2n^3}+\frac{\big(16\log n+15+16\log(2\pi)\big)^2}{2048n^4}+\frac{1}{2n^4}\bigg)\\
		=&\, 2\bigg(1-\frac{1}{2n^2}+\frac{\log n}{2n^3}+\frac{74/144+\log(2\pi)}{2n^3}+\frac{13}{4n^4}\bigg)\\
		&\,\,\cdot\bigg(1+\frac{1}{2n^2}-\frac{15/16+\log(2\pi)+\log n}{2n^3}+\frac{\big(16\log n+15+16\log(2\pi)\big)^2+1024}{2048n^4}\bigg)\,.
	\end{align*}
	It is easy to show that
	$$
	\frac{\big(16\log n+15+16\log(2\pi)\big)^2+1024}{2048n^4}\leqslant\frac{1}{32n^3}\,,
	$$
	for $n\geqslant 304$, hence
	\begin{align*}
		2 \exp\Big(&\frac{x_{n+1}-x_n}{2}\Big)\!\cosh \Big(\frac{x_{n+1}+x_n}{2}\Big)\\
		\leqslant&\,2\bigg(1-\frac{1}{2n^2}+\frac{\log n}{2n^3}+\frac{74/144+\log(2\pi)}{2n^3}+\frac{13}{4n^4}\bigg)
		\bigg(1+\frac{1}{2n^2}-\frac{14/16+\log(2\pi)+\log n}{2n^3}\bigg)\\
		\leqslant&2\bigg(1-\frac{1}{2n^2}+\frac{\log n}{2n^3}+\frac{76/144+\log(2\pi)}{2n^3}\bigg)\bigg(1+\frac{1}{2n^2}-\frac{\log n}{2n^3}-\frac{14/16+\log(2\pi)}{2n^3}\bigg),
	\end{align*}
	for $n\geqslant 468$, since then $\frac{13}{4n^4}\leqslant\frac{1}{144n^3}$. We therefore finally conclude that for $n\geqslant 396$, we have (after a straightforward computation)
	\begin{align*}
		2 \exp\Big(\frac{x_{n+1}-x_n}{2}\Big)\!\cosh \Big(\frac{x_{n+1}+x_n}{2}\Big)\!\leqslant&\,2\bigg(\!1-\frac{25}{144n^3}-\frac{1}{4n^4}+\frac{\log (n)/2+\!101/288+\log(2\pi)/2}{n^5}\bigg)\\
		&\,-2\frac{\log n+76/144+\log(2\pi)}{2n^3}\cdot\frac{\log n+14/16+\log(2\pi)}{2n^3}\\
		<&\,2\bigg(\!1-\frac{25}{144n^3}-\frac{1}{4n^4}+\frac{\log (n)/2+101/288+\log(2\pi)/2}{n^5}\bigg)\\<&\,2,
	\end{align*}
	since, by a numerical check, there holds
	$$
	-\frac{25}{144}-\frac{1}{4n}+\frac{\log (n)/2+101/288+\log(2\pi)/2}{n^2}< 0
	$$
	for each $n\geqslant 3$.
	
Since all the previous estimates are valid for $n\geqslant 468$, for such $n\in\NN$ the sequence is decreasing. By numerically checking the decreasing for $n=1,\dots,468$, we then obtain that the Lalescu sequence is {\em always} decreasing, remembering formula~\eqref{2}.

To numerically check the decreasing for $n=1,\dots,468$, we used the following code for the {\em Julia -- Version~1.10.4} programming language, with the {\em IntervalArithmetic.jl} package~\cite{IntArith}:
	
	\smallskip
	
	\begin{verbatim}
		using IntervalArithmetic
		
		function rootfactorial(n)
		I = @interval(1)
		exp = @interval(1) / @interval(n)
		for j = 1:n
		I = I * @interval(j)^exp
		end
		return I
		end
		
		function lalescu(n)
		return rootfactorial(n+1) - rootfactorial(n)
		end
		
		setdisplay(:full)
		for k =  1:500
		println("a_$(k) = $(lalescu(k))")
		if !(precedes(lalescu(k), lalescu(k-1)))
		@error("lalescu($k) is not guaranteed to be smaller 
		than lalescu($(k-1))")
		end
		end
		
	\end{verbatim}

	\begin{ackn}
		We are very grateful to Federico Poloni of the Department of Computer Science of University of Pisa for helping us with the numerical analysis part.
	\end{ackn}
	
	\bibliographystyle{amsplain}
	\bibliography{Lalbib}

\end{document}